\documentclass[12pt, a4paper]{article}
\usepackage{latexsym}
\usepackage{amsmath}
\usepackage{amssymb}
\usepackage{color}
\usepackage{amsfonts}
\usepackage{verbatim}
\usepackage{amsthm}

\usepackage{anysize,hyperref}
\input xypic
\xyoption{all}

\usepackage[perpage,symbol]{footmisc}
\renewcommand{\cal}{\mathcal}

\newcommand{\p}{\vskip 0.3cm\noindent}

\begin{document}
\title{Maximal rigid objects without loops in connected 2-CY categories are
cluster-tilting objects }
\medskip
\author{Jinde Xu \qquad Baiyu Ouyang \footnotemark[1]{}\\
{\small College of Mathematics and Computer Science},\\
{\small Key Laboratory of High Performance Computing and}\\
{\small Stochastic Information Processing (Ministry of Education of China)},\\
{\small Hunan Normal University, Changsha, Hunan 410081, P. R.
China}\\}

\date{}

\maketitle
\def\blue{\color{blue}}
\def\red{\color{red}}
\newtheorem{question}{Question}
\newtheorem{theorem}{Theorem}[section]
\newtheorem{lemma}[theorem]{Lemma}
\newtheorem{corollary}[theorem]{Corollary}

\newtheorem{proposition}[theorem]{Proposition}

\newtheorem{example}[theorem]{Example}
\newtheorem{conjecture}{Conjecture}

\theoremstyle{definition}
\newtheorem{definition}[theorem]{Definition}

\newtheorem{remark}[theorem]{Remark}
\newtheorem{remark*}[]{Remark}

\baselineskip=16pt
\parindent=0.5cm

 \footnotetext[1]{\noindent Corresponding author.\\
\hspace*{2em}Email Address: hnu\_xujinde@126.com(J.D.Xu),
 oy@hunnu.edu.cn(B.Y.Ouyang).\\
\hspace*{2em}Supported by the NSF of China {(No. 11371131)} and in
part by Construct Program of the Key Discipline in Hunan Province.}

\begin{abstract}
\baselineskip=17pt In this paper, we study the conjecture
{\textrm{II.1.9}} of \cite{birs}, which said that any maximal rigid
object without loops or 2-cycles in its quiver is a cluster-tilting
object in a connected Hom-finite triangulated 2-CY category $\cal
C$. We obtain some conditions equivalent to the conjecture, and using them we proved the conjecture.

\bigskip

\medskip
\end{abstract}

\bigskip

\section{Introduction}
\indent The theory of cluster algebras, initiated by
Fomin-Zelevinsky in \cite{fz1}, and further developed in a series of
papers, including \cite{fz2,bfz,fz3}, has turned out to have
interesting connections with many parts of algebra and other
branches of mathematics. The cluster categories associated with
finite dimensional hereditary algebras introduced in \cite{bmrrt}
and the module categories mod$\Lambda$ for $\Lambda$ the
preprojective algebra of a Dynkin quiver\cite{gls} have been
developed for the categorification of cluster algebras. This
development has both inspired new directions of investigations on
the categorical side, as well as interesting feedback on the theory
of cluster algebras. We refer the reader to the nice papers
\cite{bm,bmr1,bmr2,bmrrt,gls,keller1,reiten,ringel} for the further
results.

Both the cluster categories and the stable categories
{\underline{mod}}$\Lambda$ of preprojective algebras are
triangulated 2-Calabi-Yau  categories (2-CY for short). They both
have what are called cluster-tilting objects, which are important
since they are the analogs of clusters. Hence the theory of cluster
tilting objects in 2-CY categories has been studied in a lot of
papers, for example \cite{birs,bikr,iy,zz} etc. Cluster-tilting
objects (subcategories) in 2-CY categories have many nice
properties. {For example}, the endomorphism algebras are Gorenstein
algebras of dimension at most 1 \cite{kr}; cluster-tilting objects
have the same number of non-isomorphic indecomposable direct
summands\cite{dk}. Cluster-tilting objects are maximal rigid
objects, however the converse is not true in general. For examples
of 2-CY categories in which maximal rigid objects are not cluster
tilting please refer to \cite{bikr,bmv}. It is natural to ask under
what conditions maximal rigid objects will be cluster-tilting. In
\cite{birs}, Buan, Iyama, Reiten and Scott proved that for an exact
stably 2-CY category $\cal C$ if it admits a cluster-tilting
subcategory, then every functorially finite maximal rigid
subcategory is cluster-tilting in Theorem {$\textrm{II.1.8}$}, which was generalized to arbitrary 2-CY
categories by Zhou and Zhu (Theorem 2.6 \cite{zz}). In addition,
they gave a conjecture in \cite{birs}

\begin{conjecture}[Conjecture II.1.9 \cite{birs}]Let $\cal C$ be a connected Hom-finite triangulated $2$-CY category.
Then any maximal rigid object without loops or 2-cycles in its
quiver is a cluster-tilting object.
\end{conjecture}
The aim of this paper is to prove this conjecture.

Cluster-tilting subcategories in Calabi-Yau categories, and in
general $n$-rigid $S_n$-subcategories in triangulated categories with
a Serre functor , was systematic studied by Iyama and Yoshino
in \cite{iy}. The notion of subfactor categories was induced,
moreover Iyama and Yoshino get a result that there is a one-one
correspondence between cluster-tilting subcategories of $\cal C$
containing $\cal D$ and cluster-tilting subcategories of the
subfactor category $\cal U$ (Theorem4.9 \cite{iy}). Inspired by this
theorem, we reduce the conjecture to the following: if $T$ is an indecomposable
maximal rigid object without loops, then $T$ is a cluster-tilting object. From this
point of view, we obtain some conditions equivalent to the
conjecture, which we then use to prove
the conjecture.\p

\noindent{\textbf{Theorem 3.6'.} \textit{Let $\cal C$ be a connected
Hom-finite triangulated 2-CY category. Then any maximal rigid object
without loops in its quiver is a cluster-tilting object.} }\p

The paper is organized as follows.

In section 2, we collect some useful notions and results which will
be used in the proof of the conjecture, especially the notion of the
subfactor subcategories. In section 3, we give the proof of main
results.

\section{Preliminaries}
{Throughout this paper, $k$ denotes an algebraically closed field
and $\cal C$ denotes a $k$-linear triangulated category whose shift
functor is denoted by [1]. We assume, unless otherwise stated, that
$\cal C$ is Hom-finite and
 Krull-Schmidt, i.e. any object of $\cal C$ is
isomorphic to a finite direct sum of objects whose endomorphism
rings are local.} We denote by \textsf{ind}$\cal {C}$ the set of
isomorphism classes of indecomposable objects in $\cal C$, and by
{$\cal C(X,Y)$} the set of morphisms from $X$ to $Y$ in $\cal C$. We
denote by rad$_{\cal C}$ the \textit{Jacobson radical} of $\cal C$,
namely, rad$_{\cal C}$ is an ideal of {$\cal C$} such that
rad$_{\cal C}(X, X)$ coincides with the Jacobson radical of the
endomorphism ring End$_{\cal C}X$ for any $X\in \cal C$. For basic references
on representation theory of algebras and triangulated categories, we
refer to \cite{ass,hap}.

For a subcategory $\cal D$ of $\cal C$, we always mean that $\cal D$
is a full subcategory which is closed under isomorphisms, direct
sums and direct summands. $\cal D^{\bot}$ (resp. $^{\bot}\cal D$)
denotes the subcategory consisting of $X\in \cal C$ with $\cal
C(D,X)=0$ {(resp. $\cal C(X,D)=0$)} for any $D\in \cal D$. For any
object $T\in \cal C$ we denote by \textsf{add}$T$ the smallest
{additive} subcategory of $\cal C$ containing $T$.

A morphism $f:X\to Y$ is called \textit{right minimal} if every
$h\in$ End$_{\cal C}X$ such that $fh=f$ is an automorphism. We call $f$
\textit{a right almost split} morphism if $f\in \textrm{rad}_{\cal C}$ and
$$\xymatrix{\cal C(-,X)\ar[r]^{f_*}&rad_{\cal C}(-,Y) \ar[r]& 0}$$
is exact as functors on $\cal C$. A morphism $f$ is called
\textit{right minimal almost split} if it both right minimal and
right almost split. Dually, a \textit{left minimal almost split}
morphism is defined. For a subcategory $\cal D$ of $\cal C$, we call
$f$ \textit{a right $\cal D$-approximation} of $Y\in \cal C$ if
$X\in \cal D$ and
$$\xymatrix{\cal C(-,X)\ar[r]^{f_*}&\cal C(-,Y)\ar[r]&0}$$
is exact as functor on $\cal D$. Similarly, we call a right $\cal
D$-approximation \textit{minimal} if it is right minimal. We call
$\cal D$ \textit{a contravariantly finite subcategory} of $\cal C$
if any $Y\in \cal C$ has a right $\cal D$-approximation. Dually,
\textit{a (minimal) left $\cal D$-approximation} and \textit{a
covariantly finite subcategory} are defined. A contravariantly and
covariantly finite subcategory is called \textit{functorially
finite}. It is easy to see that \textsf{add}$T$ is functorially
finite for any object $T\in \cal C$ using Hom-finiteness and the fact that the
number of indecomposable summands of $T$ is finite.

For two subcategories $\cal X$ and $\cal Y$ of $\cal C$, $\cal {X*Y}$
denotes the collection consisting of all objects $E$ occurring in triangles $X\to E\to Y\to X[1]$, where $X\in \cal X, Y\in \cal Y$.
We call $\cal X$ extension closed if $\cal {X*X=X}$. By the
octahedral axiom, we have $\cal {(X*Y)*Z=X*(Y*Z)}$.

\begin{definition} For $X,Y\in \cal C$ and $n\in \mathbb{Z}$, we put
Ext$^n(X,Y)=\cal C(X,Y[n])$. A triangulated category $\cal C$ is
called 2-Calabi-Yau, 2-CY for short , if there are bifunctorial
isomorphisms
$$\textrm{Ext}^1(X,Y)=D\textrm{Ext}^1(Y,X)$$
for $X,Y\in \cal C$, where $D=$Hom$_k(-,k)$ is the duality of
$k$-spaces.
\end{definition}

A Hom-finite triangulated category $\cal C$ is 2-CY if and only if
it has almost split triangles with the AR-translation $\tau$ and $\tau:\cal
C\to \cal C$ is a functor isomorphic to the shift functor $[1]$ (see
\cite{rv}). An exact category is called a stably 2-CY category if it
is Frobenius, that is, it has enough projectives and injectives,
which coincide, and the stable category is a 2-CY triangulated
category. If a triangulated category is triangulated equivalent to
the stable category of a stably 2-CY exact category, then we call it an
algebraic triangulated 2-CY category{\cite{birs}}. For more examples and
information on 2-CY category please refer to
\cite{birs,keller2,keller3}.

\begin{definition}[Definition 2.1 \cite{zz}]
Let $\cal T$ be a subcategory of a triangulated 2-CY category $\cal
C$.
\begin{itemize}\setlength{\itemsep}{-0.5cm}
\item $\cal T$ is called \emph{rigid} if Ext$^1\cal {(T,T)}=0$.\\
\item $\cal T$ is called \emph{maximal rigid} if $\cal T$ is rigid and
is maximal with respect to this property, i.e. if Ext$^1(\cal T\cup
\textsf{add}M,\cal T\cup \textsf{add}M)=0$, then $M\in \cal T$.\\
\item $\cal T$ is called \emph{cluster-tilting} if $\cal T$ is
functorially finite and $\cal {T=T}[-1]^\bot=^\bot\cal T[1]$.
\end{itemize}
\end{definition}

An object $T$ is called \emph{rigid}, \emph{maximal rigid} or
\emph{cluster-tilting} if {\textsf{add}$T$} is a rigid, maximal rigid,
or cluster tilting subcategory respectively.

\begin{remark}
\begin{itemize}
\item [1.] {Higher analogous concepts of $n$-rigid, maximal $n$-rigid
and $n$-cluster-tilting subcategories were defined in \cite{iy}. As
such, rigid, maximal rigid and cluster tilting subcategories are
often known as $2$-rigid, maximal $2$-rigid and $2$-cluster-tilting
subcategories, respectively.}

\item [2.] Any triangulated 2-CY category $\cal C$ admits rigid subcategories ($0$ is viewed
as a trivial rigid object), and also admits maximal rigid
subcategories if $\cal C$ is skeletally small. But there are
triangulated 2-CY categories which contain no cluster-tilting
subcategories \cite{bikr,bmv}.
\item [3.] Cluster-tilting subcategories are functorially {finite} maximal rigid subcategories. But the converse
is not true in general \cite{bikr,bmv}. It was observed by
Buan-Marsh-Vatne \cite{bmv} that the cluster tubes contain maximal
rigid objects, but none of them are cluster-tilting objects.
\item[4.] Let $\cal T$ be a functorially finite maximal rigid subcategory of a triangulated 2-CY category $\cal C$. Then every rigid
object belongs to $\cal {T*T}[1]$. $\cal T$ is cluster-tilting if
and only if $\cal {C=T*T}[1]$ (see \cite{zz}).
\end{itemize}
\end{remark}

In order to prove our main results, let us review some useful notions
and results (for more details please refer to \cite{iy,zz}).
\begin{definition}[Definition 2.5 \cite{iy}]
Fix a subcategory $\cal D$ of $\cal C$ satisfying $\cal C(\cal
{D,D}[1])=0$. For a subcategory $\cal X$ of $\cal C$, put
$$\mu^{-1}(\cal {X;D})=(\cal {D*X}[1])\cap^\bot\cal D[1].$$
Then $\mu^{-1}(\cal {X;D})$ consists of all $C\in \cal C$ such that
there exists a triangle
$\xymatrix@C=0.3cm{X\ar[r]^f&D_X\ar[r]&C\ar[r]&X[1]}$ with $X\in \cal
X$ and a left $\cal D$-approximation $f$.\\
Dually, for a subcategory $\cal Y$ of $\cal C$, put
$$\mu(\cal {Y;D})=(\cal Y[-1]*\cal D[1])\cap\cal D[-1]^\bot.$$
Then $\mu(\cal {Y;D})$ consists of all $C\in \cal C$ such that there
exists a triangle
$\xymatrix@C=0.3cm{C\ar[r]&D_Y\ar[r]^g&Y\ar[r]&C[1]}$ with $Y\in \cal
Y$ and a right $\cal D$-approximation $g$.\\
We call a pair $(\cal {X,Y})$ of subcategories of $\cal C$ a
\emph{$\cal D$-mutation pair} if
$$\cal D\subset \cal Y \subset \mu^{-1}(\cal {X;D})\quad and\quad \cal D\subset \cal X \subset \mu(\cal {Y;D}).$$
\end{definition}

It is not difficult to see that: for subcategories $\cal {X,Y}$
containing $\cal D$, $(\cal{X,Y})$ forms a $\cal D$-mutation pair if
and only if for any $X\in \cal X$ , $Y_1\in\cal Y$ there are two
triangles:
$$\xymatrix{X\ar[r]^f&D\ar[r]^g&Y\ar[r]&X[1]}$$
$$\xymatrix{X_1\ar[r]^{f_1}&{D_1}\ar[r]^{g_1}&Y_1\ar[r]&X_1[1]}$$
where $D,D_1\in \cal D,Y\in{\cal Y},X_1\in\cal X$, $f$ and $f_1$ are
left $\cal D$-approximations; $g$ and $g_1$ are right $\cal
D$-approximations. Hence, for a $\cal D$-mutation pair $(\cal {X,Y})$,
Iyama and Yoshino construct a functor $\mathbb{G}:\cal {X/[D]\to Y/[D]}$
as follows: For any $X\in \cal X$, fix a triangle
$$\xymatrix{X\ar[r]^{\alpha_X}&D_X\ar[r]^{\beta_X}&\mathbb{G}X\ar[r]^{\gamma_X}&X[1]}$$
where $D_X\in \cal D, \mathbb{G}X\in \cal Y$, and $\alpha_X$ is a
left $\cal D$-approximation, $\beta_X$ is a right $\cal
D$-approximation, and define $\mathbb{G}X$ by this. For any morphism
$f\in {\cal C}(X,X'),\ X,X'\in {\cal X}$, there exist morphisms $g$
and $h$ which make the following diagram commutative.
$$\xymatrix{X\ar[d]^f\ar[r]^{\alpha_X}&D_X\ar[d]^g\ar[r]^{\beta_X}&\mathbb{G}X\ar[d]^h\ar[r]^{\gamma_X}&X[1]\ar[d]^{f[1]}\\
X'\ar[r]^{\alpha_{X'}}&D_{X'}\ar[r]^{\beta_{X'}}&\mathbb{G}X'\ar[r]^{\gamma_{X'}}&X'[1]}$$
Now put $\mathbb{G}\bar{f}:=\bar{h}$.

\begin{proposition}[Proposition 2.6 \cite{iy}]\label{shift} In the situation above, the following assertions hold.

\begin{itemize}
\item[\emph{(1)}] The functor $\mathbb{G}:\cal
{X/[D]\to Y/[D]}$ is an equivalence of categories.

\item[\emph{(2)}] $\cal Y =\mu^{-1}(\cal {X;D})$ and $\cal X=\mu(\cal
{Y;D})$ hold.

\end{itemize}

\end{proposition}

Next, we will introduce the notion of subfactor triangulated
categories which was also defined by Iyama and Yoshino in \cite{iy}. Let
$\cal C$ be a triangulated category and $\cal {D\subset Z}$ be
subcategories of $\cal C$. Assume $\cal Z$ and $\cal D$ satisfy the
following two conditions:
\begin{itemize}
\item[(Z1)] $\cal Z$ is extension closed, i.e. $\cal {Z*Z=Z}$.
\item[(Z2)] $\cal {(Z,Z)}$ forms a $\cal D$-mutation pair.
\end{itemize}
\begin{definition}[Definition 4.1 \cite{iy}]\label{tri}Under above setting, put the subfactor category
$$\cal U:=\cal {Z/[D]}$$
where $\cal {[D]}$ is the ideal of $\cal C$ consisting of morphisms
which factor through objects in $\cal D$.

Denote by $\langle 1\rangle$ the
equivalence $\mathbb{G}:\cal {U\to U}$ constructed in Proposition \ref{shift}.  Thus for $\forall X\in
{\cal Z}$, there exists a triangle
$$\xymatrix{X\ar[r]^{\alpha_X}&D_X\ar[r]^{\beta_X}&X\langle 1 \rangle\ar[r]^{\gamma_X}&X[1]}$$
where $\alpha_X$ is a left $\cal D$-approximation of $X$ and
$X\langle 1 \rangle\in \cal Z$. Note that $X\langle 1 \rangle$ is unique up
to summands in $\cal D$.

Let $\xymatrix{X\ar[r]^a&Y\ar[r]^b&Z\ar[r]^c&X[1]}$ be a triangle in
$\cal C$ with $X,Y,Z\in \cal Z$. Since ${\cal C}({\cal Z}[-1], {\cal
D})=0$, which makes the composition $Z[-1]\to X\to D_X$ is zero, then
there is a commutative diagram of triangles:
$$\xymatrix{X\ar[r]^a\ar@{=}[d]&Y\ar[r]^b\ar[d]&Z\ar[r]^c\ar[d]^d&X[1]\ar@{=}[d]\\
X\ar[r]^{\alpha_X}&D_X\ar[r]^{\beta_X}&X\langle 1\rangle
\ar[r]^{\gamma_X}&X[1]}$$ with $\alpha_X$ a left $\cal
D$-approximation and $\beta_X$ a right $\cal D$-approximation.
Define \emph{the triangles} in $\cal U$ to be the {diagrams} in
$\cal U$ which are isomorphic to a {diagram}
$$\xymatrix{X\ar[r]^{\overline{a}}&Y\ar[r]^{\overline{b}}&Z\ar[r]^{\overline{d}}&X\langle 1\rangle}$$
in $\cal U$, where $\overline{a},\overline{b},\overline{d}$ are the
{residue} classes of morphisms $a,b,d$.

\end{definition}

\begin{theorem}[Theorem 4.2 \cite{iy}] The category $\cal U$ forms a
triangulated category with respect to the auto-equivalence $\langle
1 \rangle$ and triangles defined in Definition \ref{tri}.

\end{theorem}

In the following, we collect some useful properties of the subfactor
category $\cal U$ which were proved in \cite{iy}.

\begin{lemma}[Proposition 4.4(1) \cite{iy}]\label{rigid}
In the situation above, if $\cal T$ is a rigid
subcategory of $\cal Z$, then so is $\overline{\cal T}$ as a
subcategory of $\cal U$.

\end{lemma}

\begin{theorem}[Theorem 4.7 of \cite{iy} in the case
$n=2$]\label{iy} Let $\cal C$ be a triangulated $2$-CY  category,
$\cal D$ a functorially finite rigid subcategory of $\cal C$. Set
$\cal Z={\cal D}[-1]^{\bot}=^{\bot} {\cal D}[1]$, then the subfactor
category $\cal {U=Z/[D]}$ forms a triangulated $2$-CY category too.

\end{theorem}

\section{Main Proof}

\begin{lemma}\label{approx}
{Let $\cal C$ be a Hom-finite Krull-Schmidt additive category,
$T=T_1\oplus T_2\oplus...\oplus T_n$ a basic object of $\cal C$,
where $T_i$ are non-isomorphic indecomposable summands of $T$.} For
any object $M \in \cal C$, if there is a non-zero minimal left
$\textsf{add}T$-approximation $f:M\to T'$ of $M$, assume
$f=\dbinom{f_1}{f''}:\xymatrix{M\ar[r]&T_1\oplus T''}$, then $f_1$ cannot factor through
$\textsf{add}(T/T_1)$.
\end{lemma}

\proof If $f_1$ factors through some non-zero object $N\in
\textsf{add}(T/T_1)$, i.e. $f_1=gh$, then we have a commutative
diagram of morphisms
$$\xymatrix{&N\oplus T''\ar[dr]^{g\ 0 \choose 0\ 1}\\
M\ar[rr]^{f_1\choose f''}\ar[ur]^{h\choose f''}& &T_1\oplus T''}$$
where $N,T'',T_1\in \textsf{add}T$.

We claim $\dbinom{h}{f''}:\xymatrix{M\ar[r]&N\oplus T''}$ is a left
$\textsf{add}T$-approximation of $M$. Indeed, For any $X\in
\textsf{add} T$, $\forall \varphi:M\to X$, there exists a morphism
$\psi:T_1\oplus T''\to X$ such that $\varphi=\psi f$ because $f$
is a left $\textsf{add}T$-approximation. Noting that
$f=\dbinom{f_1}{f''}=\dbinom{g\ \ 0}{0\ \ 1}\dbinom{h}{f''}$, hence
$\varphi=\psi\dbinom{g\ \ 0}{0\ \ 1}\dbinom{h}{f''}$. That means the
sequence
$$\xymatrix{(N\oplus T'',X)\ar[rr]^{\quad{\dbinom{h}{f''}}^*}&&(M,X)\ar[r]&0}$$
is exact.

As we know the minimal left approximation is a direct summand of any
left approximation by Proposition5.1.2 of \cite{ej}, i.e.
$T_1\oplus T''$ is a direct summand of $N\oplus T''$, $T_1$ must be a direct summand
of $N$ by the \emph{unique decomposition theorem} . But in our case, that is impossible, because $T_1$ is not
contained in $N$ as a direct summand. \qed\p

We call a subcategory $\cal T'$ an almost complete maximal rigid
subcategory if there is an indecomposable object $X$ which is not
isomorphic to any object in $\cal T'$ such that $\cal T=
\textsf{add}(\cal T'\cup \{X\})$ is a functorially finite maximal
rigid subcategory in $\cal C$. Such $X$ is called a \emph{complement} of an
almost complete maximal rigid subcategory $\cal T'$. Noting that the radical of End$_{C}X$ is nilpotent for any complement $X$ and $\cal T$ has right and left almost split morphisms by functorially finiteness of $\cal T$, it is easy to
see that any almost complete maximal rigid subcategory is
functorially finite (cf. Proposition~3.13 of \cite{as}). As before, an object $T'$ is called almost
complete maximal rigid if $\textsf{add}T'$ is an almost complete
maximal rigid subcategory.

\begin{lemma}\label{zz}\emph{[Corollary3.3 \cite{zz}]} {Let $\cal C$ be
a triangulated $2$-CY category}, $\cal T'$ an almost complete
maximal rigid subcategory of $\cal C$. Then there are exactly two
complements of $\cal T'$, say $X$ and $Y$. Denote by $\cal T=
\textsf{add}(\cal T'\cup \{X\})$, $\cal T^*= \textsf{add}(\cal
T'\cup \{Y\})$. Then $(T,T^*)$, $(T^*,T)$ are $T'$-mutation pairs.
\end{lemma}
\p

Let $T=T_1^{m_1}\oplus T_2^{m_2}\oplus...\oplus T_n^{m_n}$ be an
object in $\cal C$, where the $T_i$ are pairwise non-isomorphic
indecomposable objects and $m_i\geq 1$. Let $S_i=S_{T_i}$ be the
simple End$_\cal C(T)$-module corresponding to $T_i$, and $P_i=\cal
C(T_i,T)$ be the indecomposable projective End$_\cal C(T)$-module
with top $S_i$. {It is well known that} the following numbers are
equal for $1\leq i,j\leq n$ {\cite{ass,gls}}:
\begin{itemize}
\item The number of arrows $i\to j$ in the quiver of End$_\cal CT$;
\item dim Ext$^1_{End(T)}(S_i,S_j)$;
\item The dimension of the space of irreducible maps $T_i\to T_j$ in the
category \textsf{add}$T$
\end{itemize}

\begin{lemma}\label{loop}
Let $X\ncong Y$ be the two complements of a basic almost complete
maximal rigid object $T'$. Then the following are equivalent:
\begin{itemize}
\item The quiver of End$_\cal C(T'\oplus X)$ has no loops at $X$;
\item Every non-isomorphism $\varphi: X\to X$ factors through $\emph{\textsf{add}}T'$;
\item dim$\,Ext^1_\cal C(Y,X)=1$.
\end{itemize}
\end{lemma}

\proof The proof of Lemma 6.1 in \cite{gls} works also in this setting. For the
convenience of the reader we briefly give the proof. The equivalence
of the first two statements is easy to show. By Lemma \ref{zz},
there is a triangle
$$\xymatrix@C=1.2cm{X\ar[r]^f&T_1\ar[r]^g&Y\ar[r]&X[1]}$$
where $f$ is a left approximation and $T_1\in \textsf{add} T'$.
Applying ${\cal C}(-,X)$ functor yields an exact sequence
$$\xymatrix@C=1.2cm{{\cal C}(T_1,X)\ar[r]^{{\cal C}(f,X)}&{\cal C}(X,X)\ar[r]&{\cal C}(Y[-1],X)\ar[r]&0}$$
Since $f$ is an $\textsf{add} T'$-approximation, every
non-isomorphism $\varphi: X\to X$ factors through $\textsf{add} T'$
if and only if it factors through $f$. This is equivalent to the
cokernel $Ext^1_\cal C(Y,X)=1$ of ${\cal C}(f,X)$ being
1-dimensional. Here we use that $k$ is an algebraically closed field, which  implies that $\cal C(X,X)/\textrm{rad}_{\cal C}(X,X)\cong k$. \qed

\begin{lemma}\label{ar}
Let $\cal C$ be a triangulated 2-CY category, $T$ a basic maximal
rigid object without loops in its quiver. Assume $T=T_1\oplus
T_2\oplus...\oplus T_n$ where $T_i$ are non-isomorphic
indecomposable summands of $T$ and {$T'=T/T_i$} is an
almost complete maximal rigid object {for any $1\leq i\leq n$}.\\
Put $\cal D=\emph{\textsf{add}}T'\subset\cal Z=\cal D[-1]^\bot=^\bot
\cal D[1],$ and $\cal U=\cal Z/[\cal D]$. Then
$$\xymatrix{T\langle 1\rangle\ar[r]&0\ar[r]&T\ar[r]&T}\qquad and \qquad \xymatrix{T\ar[r]&0\ar[r]&T\langle 1\rangle\ar[r]&T\langle 1\rangle}$$
are AR-triangles in $\cal U$.
\end{lemma}

\proof By {Theorem \ref{iy}}, $\cal U$ is a 2-CY category with shift
$\langle 1\rangle$, and $0$ is an almost complete maximal rigid
object. By Lemma \ref{zz}, there are exactly two complements of $0$,
that is to say there are only two maximal rigid objects (up to
isomorphism) in $\cal U$, which consist of exactly one direct
summand. Because $\cal U(T\langle 1\rangle,T\langle 2\rangle)=\cal
U(T,T\langle 1\rangle)=0$ by Proposition \ref{rigid}, $T$ and
$T\langle 1\rangle$ are maximal rigid objects of $\cal U$, similarly
$T\langle i\rangle$ is a maximal rigid object for any $i\in
\mathbb{Z}$. But there are only two maximal rigid objects in $\cal
U$ and $T\langle i\rangle$ and $T\langle i+1\rangle$ are distinct
maximal rigid objects by $\cal U(T\langle i\rangle,T\langle
i+1\rangle)=0$, hence we get that
$$T\langle -1\rangle \cong T\langle 1\rangle\cong ...\cong T\langle
odd\rangle\qquad and \qquad T\langle 0\rangle \cong T\langle
2\rangle\cong ...\cong T\langle even\rangle$$ Because the quiver of
End$_\cal C T$ has no loops, clearly the quiver of End$_\cal U T$ has
no loops too. By Lemma \ref{loop}, dim\,Ext$^1_\cal U(T,T\langle
1\rangle)=$dim\,Ext$^1_\cal U(T\langle 1\rangle,T)=1$, Noting that
$\langle 1\rangle=\tau_\cal U$ because $\cal U$ is 2-CY, then
dim\,Ext$^1_\cal U(T,\tau T)=$dim\,Ext$^1_\cal U(\tau T,T)=1$. That
means the non-split triangles
$$\xymatrix{T\langle 1\rangle\ar[r]&0\ar[r]&T\ar[r]&T}\qquad and \qquad \xymatrix{T\ar[r]&0\ar[r]&T\langle 1\rangle\ar[r]&T\langle 1\rangle}$$
are AR-triangles in $\cal U$. \qed

\begin{theorem}\label{tilting}
{Let $\cal C$ be a Hom-finite triangulated 2-CY category}, $T$ a
basic maximal rigid object whose quiver has no loops,
$\overline{\cal D}=T[-1]^\bot=^\bot T[1]$. Then the following are
equivalent.
\begin{itemize}
\item[]\verb"(1)" $T$ is a cluster-tilting object;
\item[]\verb"(2)" $^\bot T\cap \overline{\cal D}\cap T^\bot=0$;
\item[]\verb"(3)" $^\bot T\cap \overline{\cal D}=0$ or $\overline{\cal D}\cap
T^\bot=0$.
\end{itemize}

\end{theorem}

\proof \begin{itemize}
\item[] Obviously $(3)\Rightarrow (2)$.

\item[] $(1)\Rightarrow (3)$: By the definition of cluster-tilting, $\overline{\cal
D}=\textsf{add}T$, hence $$^\bot T\cap \overline{\cal D}=^\bot T\cap
\textsf{add}T=0.$$ Similarly $\overline{\cal D}\cap T^\bot=0$.

\item[] $(2)\Rightarrow (1)$: In order to prove $T$ is cluster
tilting, it {suffices} to prove $\overline{\cal D}\subseteq
\textsf{add}T$.

Let $X$ be a non-zero indecomposable object of $\overline{\cal D}$.
Then $\cal C(X,T)\neq0$ or $\cal C(T,X)\neq 0$ by condition (2).

\ \  If ${{\cal C}}(X,T)\neq 0$, then we have a non-zero minimal
left $\textsf{add} T$-approximation $f:X\to T'$ of X. Because $f\neq
0$, by decomposing $T'$ we can find an indecomposable summand $T_i\neq
0$ and $f_i\neq 0$ such that
$$f=\dbinom{f_i}{f''}:\xymatrix{X\ar[r]&T_i\oplus T''}.$$
We can assume $i=1$, then
$f=\dbinom{f_1}{f''}:\xymatrix{X\ar[r]&T_1\oplus T''}$, and
$f_1:X\to T_1$ does not factor through \textsf{add}$(T/T_1)$ by Lemma
\ref{approx}.

Now put $\cal D=\textsf{add}(T/T_1),\ \cal D\subset\cal Z=\cal
D[-1]^\bot=^\bot \cal D[1]$ and $\cal {U=Z/[D]}$. Noting $f_1\in
\cal Z$ because $\overline{\cal D}\subseteq \cal Z$, then the
{residue} of $f_1$ in $\cal U$, denoted by $\bar{f_1}$, is not zero
by the definition of $\cal U$ ($f_1$ does not factor through
\textsf{add}$(T/T_1)$), and it follows that $X\neq 0$ in $\cal U$.\\
Then we get that $X\in \textsf{add}T$ in $\cal C$ as follows. By
Lemma \ref{ar}, the triangle
$$\xymatrix{T_1\langle
1\rangle\ar[r]&0\ar[r]&T_1\ar[r]&T_1}$$ is the AR-triangle ending
{at} $T_1$, noting $T_1\cong T$ in $\cal U$. {If $\bar{f_1}$ were
not a retraction, then it would factor though $0\to T_1$ (right
minimal almost split)}, i.e. there is a $g$ such that
${\bar{f}_1}=0\cdot g=0$, in contradiction to $\bar{f_1}\neq 0$.
Therefore, $\bar{f_1}$ is a retraction, but $X$ is
indecomposable, then $X\cong T_1\cong T$ in the subfactor category
$\cal U$. {That means $X\oplus D_1\cong T_1\oplus D_2$  in $\cal Z$,
$D_1,D_2\in \cal D$, by the \emph{unique decomposition theorem} and
$T_1\notin \cal D$, $T_1$ must be a direct summand of $X$. Again
because $X$ is indecomposable in $\cal Z$, we get $X\cong T_1$ in
$\cal Z$ }, hence in $\cal C$, that means $X\in \textsf{add} T$.

\ \  If $\cal C(T,X)\neq 0$, by the same argument using the dual of Lemma
\ref{approx}, we also have $X\in \textsf{add} T$. Hence,
$\overline{\cal D}\subseteq \textsf{add} T$, i.e. $T$ is a cluster-tilting object. We have completed the proof. \qed

\end{itemize}
\p

Now, we can prove the conjecture. Recall that a triangulated category
$\cal C$ is said to be connected if the underlying graph
$\overline{\Gamma}_{\cal C}$ of the AR-quiver $\Gamma_{\cal C}$ is a
connected graph, where \emph{the underlying graph
$\overline{\Gamma}_{\cal C}$} is obtained from $\Gamma_{\cal C}$ by
forgetting the orientation of the arrows. For the definition of AR-quiver of a triangulated category please refer to \cite{hap}.
\begin{theorem}\label{main}
Let $\cal C$ be a connected Hom-finite triangulated 2-CY category.
Then any maximal rigid object $T$ without loops or 2-cycles in its
quiver is a cluster-tilting object.
\end{theorem}

\proof Let $T$ be a basic maximal rigid object whose quiver has no
loops. Using the notations from Theorem {\ref{tilting}}, it suffices
to show that $^\bot T\cap \overline{{\cal D}}\cap T^\bot=0$.

If $^\bot T\cap \overline{{\cal D}}\cap T^\bot\neq 0$, let $0\neq
M\in ^\bot T\cap \overline{{\cal D}}\cap T^\bot$ with $M$
indecomposable. If the space of irreducible maps Irr$_{\cal C}(N,M)\neq 0$ for some indecomposable object
{$N$}, i.e. there exists an irreducible map from $N$ to $M$. We
claim that $N$ is also in $^\bot T\cap \overline{{\cal D}}\cap
T^\bot$. Indeed, taking the AR-triangle ending at $M$
$$\xymatrix{\tau M\ar[r]&E\ar[r]&M\ar[r]&\tau M[1]},$$
and applying functors $\cal C(-,T[1])$ and $\cal C(T,-)$, we obtain
the following exact sequences
\begin{gather}
\xymatrix{0=\cal C(M,T[1])\ar[r]&\cal C(E,T[1])\ar[r]&\cal C(\tau
M,T[1])=\cal C(M[1],T[1])=0}\tag{a}
\end{gather}

\begin{gather}
\xymatrix{0=\cal C(T,M[1])=\cal C(T,\tau M)\ar[r]&\cal
C(T,E)\ar[r]&\cal C(T,M)=0} \tag{b}
\end{gather}

Then we get $E\in \overline{{\cal D}}$ from (a) and $E\in T^\bot$
from (b).

To prove our claim, it suffices to show $E\in ^\bot T$. If ${\cal
C}(E,T)\neq 0$, there at least exists a non-zero indecomposable
direct summand $X$ of $E$ such that ${\cal C}(X,T)\neq 0$. {Noting
that $X\in \overline{{\cal D}}$ since $E\in \overline{{\cal D}}$,
then using the same arguments in the proof of $``(2)\Rightarrow (1)"$ in
Theorem \ref{tilting}, we get $X\in \textsf{add} T$ in $\cal C$.
Since $X\in \textsf{add} T$ is an indecomposable direct summand of
$E$, there exist irreducible morphisms between $X$ and $M$. Then
${\cal C}(T, M)\neq 0$, in contradiction to $M\in T^\bot$.} Hence $E\in ^{\bot}T$ and thus, $N\in
^\bot T\cap \overline{{\cal D}}\cap T^\bot$.

If Irr$_{\cal C}(M,N')\neq 0$, dually we can prove $N'\in^\bot T\cap
\overline{{\cal D}}\cap T^\bot$. (Applying functors ${\cal
C}(T[-1],-)$ and ${\cal C}(-,T)$ to the AR-triangle $\xymatrix{
M\ar[r]&E'\ar[r]&M[-1]\ar[r]&M[1]}.$)

We have shown that if $N$ (resp. $N'$) is a direct
predecessor (resp. successor) of $M\in ^\bot T\cap \overline{{\cal D}}\cap
T^\bot$, then $N$ (resp. $N'$) also belongs to $^\bot T\cap \overline{{\cal
D}}\cap T^\bot.$ That means the connected component which contains
$M$ is contained in $^\bot T\cap \overline{{\cal D}}\cap T^\bot$, which is of
course contained in $^\bot T$. So $T$ and $M$ cannot be in the same
connected component, in contradiction to the connectness of $\cal C$.
\qed \p

Noting that in the proof of the theorem, we do not use the
condition that the quiver of End$_{\cal C}T$ has no 2-cycles. Hence the
Theorem \ref{main} can be restated as follows:\p

\noindent \textbf{Theorem 3.6'.} \textit{Let $\cal C$ be a connected
Hom-finite triangulated 2-CY category. Then any maximal rigid object
without loops in its quiver is a cluster-tilting object.} \p

The following result is an immediate consequence of Theorem 3.6' by
combining Theorem 2.6 of [ZZ]

\begin{corollary}
Let $\cal C$ be a connected Hom-finite triangulated 2-CY category. If $\cal C$ has
a maximal rigid object whose quiver has no loops, then every maximal rigid
object is cluster-tilting.
\end{corollary}

\begin{example}
Let $R$ be a one-dimension simple hypersurface singularity. Then the
category $CM(R)$ of maximal Cohen-Macaulay modules is a Frobenius
category, and so that the stable category $\underline{CM}(R)$ is a
triangulated 2-CY category \cite{bikr}. In the case $D_n$ with $n$
odd, the AR-quiver of $\underline{CM}(R)$ is the following
(\cite{bikr} or \cite {y})

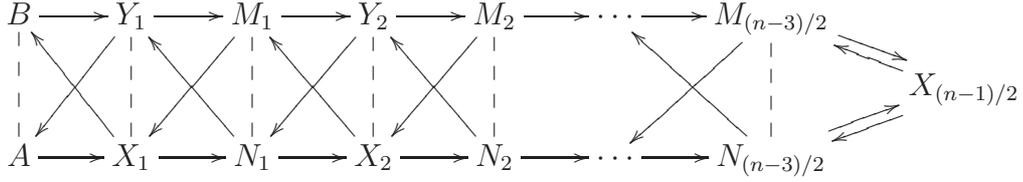
$$\xymatrix@R=0.3cm{B\ar@{--}[dd]\ar[r]&Y_1\ar[ldd]\ar@{--}[dd]\ar[r]&M_1\ar[ldd]\ar@{--}[dd]\ar[r]&Y_2\ar[ldd]\ar@{--}[dd]\ar[r]
&M_2\ar[ldd]\ar@{--}[dd]\ar[r]&\cdots\ar[r]&M_{(n-3)/2}\ar[ldd]\ar@{--}[dd]\ar@<0.4ex>[rd]\\
&&&&&&&X_{(n-1)/2}\ar@<0.4ex>[lu]\ar@<0.4ex>[ld]\\
A\ar[r]&X_1\ar[luu]\ar[r]&N_1\ar[luu]\ar[r]&X_2\ar[luu]\ar[r]&N_2\ar[luu]\ar[r]&\cdots\ar[r]&N_{(n-3)/2}\ar[luu]\ar@<0.4ex>[ru]}$$
where the dotted line between two indecomposable modules means that
they are connected via $\tau$. Proposition 2.5 of \cite{bikr} proved
that the indecomposable modules $A$ and $B$ are maximal rigid
objects in $\underline{CM}(R)$ but none of them is cluster-tilting.
And the computation of the support of $\underline{Hom}(A,-)$ is the
following
$$
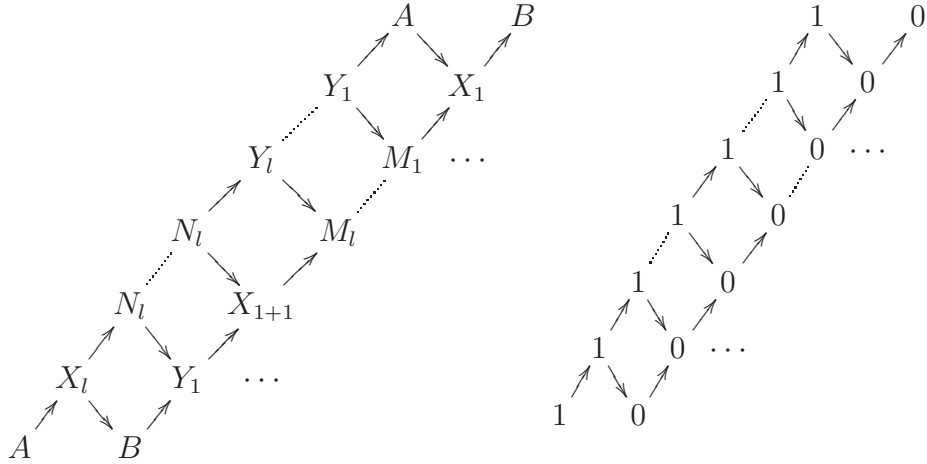
\xymatrix@R=0.4cm@C=0.1cm{&&&&&&A\ar[dr]&&B\\
&&&&&Y_1\ar[ur]\ar[dr]&&X_1\ar[ur]\\
&&&&Y_l\ar[dr]\ar@{.}[ur]&&M_1\ar[ur]&\cdots\\
&&&N_l\ar[dr]\ar[ur]&&M_l\ar@{.}[ur]\\
&&N_l\ar[dr]\ar@{.}[ur]&&X_{1+1}\ar[ur]\\
&X_l\ar[dr]\ar[ur]&&Y_1\ar[ur]&\cdots\\
A\ar[ur]&&B\ar[ur]} \xymatrix@R=0.4cm@C=0.1cm{&&&&&&1\ar[dr]&&0\\
&&&&&1\ar[ur]\ar[dr]&&0\ar[ur]\\
&&&&1\ar[dr]\ar@{.}[ur]&&0\ar[ur]&\cdots\\
&&&1\ar[dr]\ar[ur]&&0\ar@{.}[ur]\\
&&1\ar[dr]\ar@{.}[ur]&&0\ar[ur]\\
&1\ar[dr]\ar[ur]&&0\ar[ur]&\cdots\\
1\ar[ur]&&0\ar[ur]}$$ where $B=\tau A$ and $l=(n-3)/2$. Hence the
quiver of $\underline{End}(A)$ has a loop, the case of $B$ is the
same.
\end{example}

But the condition that the quiver has no loops is not necessary for
a cluster-tilting object, for example:

\begin{example}
Let $R$ be a one-dimension simple hypersurface singularity in the
case $A_n$ with $n$ odd. the AR-quiver of $\underline{CM}(R)$ is the
following (\cite{bikr} or \cite {y})
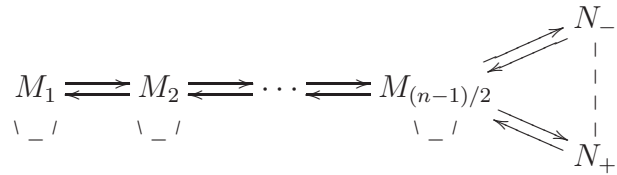
$$\xymatrix@R=0.3cm{&&&&N_{-}\ar@<0.4ex>[dl]\ar@{--}[dd]\\
M_1\ar@(dl,dr)@{--}\ar@<0.4ex>[r]&M_2\ar@(dl,dr)@{--}\ar@<0.4ex>[l]\ar@<0.4ex>[r]&\ar@<0.4ex>[l]\cdots\ar@<0.4ex>[r]&M_{(n-1)/2}\ar@(dl,dr)@{--}\ar@<0.4ex>[l]\ar@<0.4ex>[ur]\ar@<0.4ex>[dr]\\
&&&&N_{+}\ar@<0.4ex>[ul]}$$ It was proved in Proposition 2.4 of
\cite{bikr} that $N_{+}$ and $N_{-}$ are cluster-tilting objects in
$\underline{CM}(R)$, and the computation of the support of
$\underline{Hom}(N_{-},-)$ is the following
$$
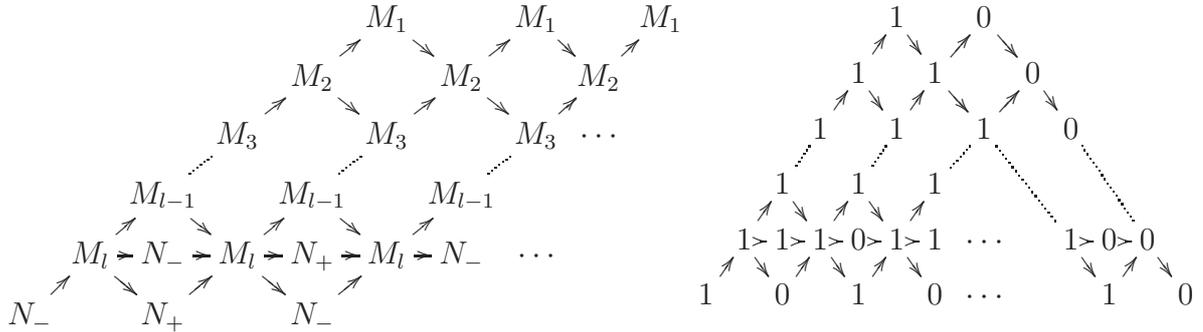
\xymatrix@R=0.5em@C=0.1em{&&&&&M_1\ar[dr]&&M_1\ar[dr]&&M_1\\
&&&&M_2\ar[dr]\ar[ur]&&M_2\ar[dr]\ar[ur]&&M_2\ar[ur]\\
&&&M_3\ar[ur]&&M_3\ar[ur]&&M_3\ar[ur]&\cdots\\
&&M_{l-1}\ar[dr]\ar@{.}[ur]&&M_{l-1}\ar[dr]\ar@{.}[ur]&&M_{l-1}\ar@{.}[ur]\\
&M_l\ar[dr]\ar[r]\ar[ur]&N_{-}\ar[r]&M_l\ar[r]\ar[dr]\ar[ur]&N_{+}\ar[r]&M_l\ar[ur]\ar[r]&N_{-}&\cdots\\
N_{-}\ar[ur]&&N_{+}\ar[ur]&&N_{-}\ar[ur]} \xymatrix@R=0.6em@C=0.2em{&&&&&1\ar[dr]&&0\ar[dr]\\
&&&&1\ar[dr]\ar[ur]&&1\ar[dr]\ar[ur]&&0\ar[dr]\\
&&&1\ar[ur]&&1\ar[ur]&&1\ar@{.}[ddrr]\ar[ur]&&0\ar@{.}[ddrr]\\
&&1\ar[dr]\ar@{.}[ur]&&1\ar[dr]\ar@{.}[ur]&&1\ar@{.}[ur]\\
&1\ar[dr]\ar[r]\ar[ur]&1\ar[r]&1\ar[r]\ar[dr]\ar[ur]&0\ar[r]&1\ar[dr]\ar[ur]\ar[r]&1&\cdots&&1\ar[dr]\ar[r]&0\ar[r]&0\ar[dr]\\
1\ar[ur]&&0\ar[ur]&&1\ar[ur]&&0&\cdots&&&1\ar[ur]&&0}$$ where
$l=(n-1)/2$. This shows that the quiver of End$(N_{-})$ has a loop
even though $N_{-}$ is cluster-tilting.

\end{example}

 \section*{ACKNOWLEDGMENTS}

The authors would like to thank Professor Bin Zhu for his helpful
suggestions in this work. The authors also would like to thank the
referee for the very kind and helpful comments and advice in shaping
the paper into its present form.

\end{document}